\documentclass[12 pt a4 paper twoside]{article}
\usepackage{amsmath,amssymb}
\usepackage{longtable,mathrsfs}
\usepackage{hyperref}
\numberwithin{equation}{section}
\newtheorem{theorem}{Theorem}[section]
\newtheorem{definition}{Definition}[section]
\newtheorem{remark}{Remark}[section]
\newtheorem{lem}{Lemma}[section]

\numberwithin{equation}{section}

\newcommand{\R}{\mathbb{R}}
\newcommand{\N}{\mathbb{N}}
\allowdisplaybreaks
\begin{document}
\begin{center}
{\textbf{Existence and uniqueness of solution of Cauchy-type problem for Hilfer fractional differential equations}}\\
\quad

    \textbf{D. B. Dhaigude and $^{*}$Sandeep P. Bhairat}\\
    Department of Mathematics,\\ Dr. Babasaheb Ambedkar Marathwada University,\\ Aurangabad -- 431 004, (M.S) India.\\
    dnyanraja@gmail.com and $^{*}$sandeeppb7@gmail.com
\end{center}
\begin{abstract}
\noindent The Cauchy-type problem for a nonlinear differential equation involving Hilfer fractional derivative is considered. We prove existence, uniqueness and continuous dependence of a solution for Cauchy-type problem using successive approximations and generalized Gronwall inequality.
\end{abstract}
\footnote{${^{*}}$ Author for correspondence, email: sandeeppb7@gmail.com.}
\section{INTRODUCTION}
Recently, the fractional differential equations have received much attention as they have emerged in various applications of science and engineering from the study of exact description of nonlinear phenomena. It has been found that the models using mathematical tools from fractional calculus describe various \break{complex} phenomena such as control, viscoelasticity, dielectric relaxations, electrochemistry, bioengineering, porous media, and many other branches of science. A huge amount of mathematically and physically interesting work, including several excellent monographs, has been published in the literature since the end of $20^{th}$ century, see \cite{hr,kst,fm,pi,skm} and references therein.

The investigation of basic theory of fractional differential equations involves the existence and uniqueness of solutions on the finite interval. The numerous results on existence and uniqueness of solution of fractional differential equations under different conditions are studied in \cite{dn,kmf,hlt,kbt,km,km1,cj,lv,jd1,jd2}, by different techniques. Recently, in \cite{kmf}, Furati et.al. obtained the existence and uniqueness of the initial value problem (IVP) for the fractional differential equation
\begin{equation}\label{s1}
D_{a^+}^{\alpha,\beta}y(x)=f(x,y),\qquad 0<\alpha<1,\, 0\leq\beta\leq1,
\end{equation}
\begin{equation}\label{s2}
  I_{a^+}^{1-\gamma}y(a)=y_a,\qquad \gamma=\alpha+\beta(1-\alpha),
\end{equation}
by using Banach fixed point technique, where $D_{a^+}^{\alpha,\beta}$ is the Hilfer (generalized Riemann-Liouville) fractional derivative of order $\alpha$ and type $\beta,$ see \cite{rh}. The fixed point technique does not indicate the interval of existence of solution, which is a necessary aspect for application purposes. This drawback of fixed point technique is removed by using Picard's iterative technique and existence of a solution is confirmed on $[a,b].$

In this paper, we study the existence, uniqueness and continuous dependence of general Cauchy-type problem by Picard's successive approximations. Clearly, the problem is well posed and investigates the qualitative properties of solution of Cauchy-type problem associated to Hilfer fractional differential equations.

The rest of the paper is organized as follows: in section 2, some preliminary results and notations are provided. The existence and uniqueness results are proved in section 3. The dependence of solutions on order and initial conditions is studied in the last section.
\section{Preliminaries}
In this section, we present some definitions, lemmas and weighted spaces which are useful in further development of this paper. For more details, see \cite{kst}.

Let $-\infty<a<b<+\infty.$ Let $C[a,b], AC[a,b]$ and $C^{n}[a,b]$ be the spaces of continuous, absolutely continuous, $n-$times continuous and continuously differentiable functions on $[a,b],$ respectively. Here $L^{p}(a,b), p\geq1,$ is the space of Lebesgue integrable functions on $(a,b).$ Further more we recall following weighted spaces \cite{kst}:
\begin{align}\label{e1}
  C_{\gamma}[a,b]&=\{f:(a,b]\to\R :  (x-a)^{\gamma}f(x)\in C[a,b]\},\quad0\leq\gamma<1,\nonumber\\
  C_{1-\gamma}[a,b]&=\{f:(a,b]\to\R : (x-a)^{1-\gamma}f(x)\in C[a,b]\},\quad0\leq\gamma<1,\\
  C_{\gamma}^{n}[a,b]&=\{f:(a,b]\to\R, f\in C^{n-1}[a,b]:f^{(n)}(x)\in C_{\gamma}[a,b]\},\, n\in \N.\nonumber
\end{align}
Clearly, $D_{a^+}^{\alpha,\beta}f=I_{a^+}^{\beta(1-\alpha)}D_{a^+}^{\gamma}f$ and $C_{1-\gamma}^{\gamma}[a,b]\subset C_{1-\gamma}^{\alpha,\beta}[a,b],\,\gamma=\alpha+\beta-\alpha\beta$, $0<\alpha<1, 0\leq\beta\leq1,$ see \cite{kmf}.
Consider the space $ C_{\gamma}^{0}[a,b]$ with the norm
\begin{equation}\label{n1}
{\|f\|}_{C_{\gamma}^{n}}=\sum_{k=0}^{n-1}{\|f^{(k)}\|}_{C}+{\|f^{(n)}\|}_{C_{\gamma}},
\end{equation}
and $C_{1-\gamma}[a,b]$ is complete metric space with the metric $d$ defined by\\
$d(y_1,y_2)={\|y_1-y_2\|}_{C_{1-\gamma}[a,b]}:=\max_{x\in[a,b]}|(x-a)^{1-\gamma}[y_1(x)-y_2(x)]|,$\\ for details see \cite{kmf}.
\begin{definition} \cite{cj}
Let $\Omega=(a,b]$ and $f:(0,\infty)\to\R$ is a real valued continuous function. The Riemann-Liouville fractional integral of a function $f$ of order $\alpha\in{\R}^{+}$ is denoted as $I_{a^+}^{\alpha}f$ and defined by
\begin{equation}\label{d1}
I_{a^+}^{\alpha}f(x)=\frac{1}{\Gamma(\alpha)}\int_{a}^{x}\frac{f(t)dt}{(x-t)^{1-\alpha}},\quad x>a,
\end{equation}
where $\Gamma(\alpha)$ is the Euler's Gamma function.
\end{definition}
\begin{definition} \cite{kst}
Let $\Omega=(a,b]$ and $f:(0,\infty)\to\R$ is a real valued continuous function. The Riemann-Liouville fractional derivative of function $f$ of order $\alpha\in{\R}_{0}^{+}=[0,+\infty)$ is denoted as $D_{a^+}^{\alpha}f$ and defined by
\begin{equation}\label{d2}
D_{a^+}^{\alpha}f(x)=\frac{1}{\Gamma(n-\alpha)}\frac{d^{n}}{dx^{n}}\int_{a}^{x}\frac{f(t)dt}{(x-t)^{\alpha-n+1}},
\end{equation}
where $n=[\alpha]+1,$ and $[\alpha]$ means the integral part of $\alpha,$ provided the right hand side is pointwise defined on $(0,\infty).$
\end{definition}
\begin{definition} \cite{hr} The Hilfer fractional derivative $D_{a^+}^{\alpha,\beta}$ of function $f\in L^{1}(a,b)$ of order $n-1<\alpha<n$ and type $0\leq\beta\leq1$ is defined by
\begin{equation}\label{d3}
    D_{a^+}^{\alpha,\beta}f(x)=I_{a^+}^{\beta(n-\alpha)}D^{n}I_{a^+}^{(1-\beta)(n-\alpha)}f(x),
\end{equation}
where $I_{a^+}^{\alpha}$ and $D_{a^+}^{\alpha}$ are Riemann-Liouville fractional integral and derivative defined by \eqref{d1} and \eqref{d2}, respectively.
\end{definition}
\begin{definition}\cite{cj}
Assume that $f(x,y)$ is defined on set $(a,b]\times G, G\subset\R.$ A function $f(x,y)$ satisfies Lipschitz condition with respect to $y,$ if for all $x\in(a,b]$ and for $y_1,y_2\in G,$
\begin{equation}\label{d4}
  |f(x,y_1)-f(x,y_2)|\leq A|y_1-y_2|,
\end{equation}
where $A>0$ is Lipschitz constant.
\end{definition}
\begin{definition}\cite{kmf} Let $0<\alpha<1,0\leq\beta\leq1,$ the weighted space $C_{1-\gamma}^{\alpha,\beta}[a,b]$ is defined by
\begin{equation}\label{w1}
C_{1-\gamma}^{\alpha,\beta}[a,b]= \big\{ f\in {C_{1-\gamma}[a,b]}:D_{a^+}^{\alpha,\beta}f \in {C_{1-\gamma}[a,b]} \big\}, \quad\gamma=\alpha+\beta(1-\alpha).
\end{equation}
\end{definition}
\begin{lem} \cite{kmf}
If $\alpha>0$ and $0\leq\mu<1,$ then $I_{a^{+}}^{\alpha}$ is bounded from $C_{\mu}[a,b]$ into $C_{\mu}[a,b].$
In addition, if $\mu\leq\alpha,$ then $I_{a^{+}}^{\alpha}$ is bounded from $C_{\mu}[a,b]$ into $C[a,b].$
\end{lem}
\begin{lem} \cite{kst}
For $x>a,$ we have
\begin{description}
\item[(i)] $I_{a^{+}}^{\alpha}(x-a)^{\beta-1}=\frac{\Gamma(\beta)}{\Gamma(\beta+\alpha)}{(x-a)}^{\beta+\alpha-1}, \quad\alpha\geq0,\beta>0.$\\
\item[(ii)] $ D_{a^{+}}^{\alpha}(x-a)^{\alpha-1}=0,\quad\alpha\in(0,1).$
\end{description}
\end{lem}
\begin{lem}\cite{kst}
Let $0<\alpha<1,0\leq\mu<1.$ If $f\in C_{\mu}[a,b]$ and $I_{a^{+}}^{1-\alpha}f\in C_{\mu}^{1}[a,b],$ then
\begin{equation*}
I_{a^{+}}^{\alpha}D_{a^{+}}^{\alpha}f(x)=f(x)-\frac{I_{a^{+}}^{1-\alpha}f(a)}{\Gamma(\alpha)}(x-a)^{\alpha-1},  \qquad\text{for all}\quad x\in(a,b].
\end{equation*}
\end{lem}
\begin{lem}\label{le}\cite{kmf} Let $\gamma=\alpha+\beta-\alpha\beta$ where $0<\alpha<1$ and $0\leq\beta\leq1.$ Let $f:(a,b]\times\R\to\R$ be a function such that $f(x,y)\in C_{1-\gamma}[a,b]$ for any $y\in C_{1-\gamma}[a,b].$ If $y\in C_{1-\gamma}^{\gamma}[a,b],$ then $y$ satisfies IVP \eqref{s1}-\eqref{s2} if and only if $y$ satisfies the Volterra integral equation of second kind
\begin{equation}\label{s3}
y(x)=\frac{y_a}{\Gamma(\gamma)}(x-a)^{\gamma-1}+\frac{1}{\Gamma(\alpha)}\int_{a}^{x}(x-t)^{\alpha-1}f(t,y(t))dt,\quad x>a.
\end{equation}
\end{lem}
\section{Existence and Uniqueness}
In this section we prove the existence and uniqueness of solution of Cauchy-type problem \eqref{s1}-\eqref{s2} in $C_{1-\gamma}^{\alpha,\beta}[a,b].$ We need the following lemma.
\begin{lem} If $\gamma=\alpha+\beta-\alpha\beta,0<\alpha<1,0\leq\beta\leq1,$ then the Riemann-Liouville fractional integral operator $I_{a^{+}}^{\alpha}$ is bounded from $C_{1-\gamma}[a,b]$ to $C_{1-\gamma}[a,b]:$
\begin{equation}\label{l1}
{\|I_{a^{+}}^{\alpha}f\|}_{C_{1-\gamma}[a,b]}\leq M\frac{\Gamma(\gamma)}{\Gamma(\gamma+\alpha)}(x-a)^{\alpha},
\end{equation}
where, $M$ is the bound of a bounded function $f.$
\end{lem}
\textbf{Proof:} From Lemma 2.1, the result follows. Now we prove the estimate \eqref{l1}. By the weighted space given in \eqref{e1}, we have
\begin{align*}
 {\|I_{a^{+}}^{\alpha}f\|}_{C_{1-\gamma}[a,b]} &= {\|(x-a)^{1-\gamma}I_{a^{+}}^{\alpha}f\|}_{C[a,b]}\\
  & \leq {\|f\|}_{C_{1-\gamma}[a,b]}{\|I_{a^{+}}^{\alpha}(x-a)^{\gamma-1}\|}_{C_{1-\gamma}[a,b]},
 \end{align*}
using Lemma 2.2, we get
 \begin{equation*}
 {\|I_{a^{+}}^{\alpha}f\|}_{C_{1-\gamma}[a,b]}\leq M \frac{\Gamma(\gamma)}{\Gamma(\alpha+\gamma)}(x-a)^{\alpha}.
\end{equation*}
Thus the proof is complete.
\begin{theorem} Let $\gamma=\alpha+\beta-\alpha\beta$ where $0<\alpha<1$ and $0\leq\beta\leq1.$ Let $f:(a,b]\times\R\to\R$ be a function such that $f(x,y)\in C_{1-\gamma}[a,b]$ for any $y\in C_{1-\gamma}[a,b],$ and satisfies Lipschitz condition \eqref{d4} with respect to $y.$ Then there exists a unique solution $y(x)$ for the Cauchy-type problem \eqref{s1}-\eqref{s2} in $C_{1-\gamma}^{\alpha,\beta}[a,b].$
\end{theorem}
\textbf{Proof:} The integral equation \eqref{s3} makes sense in any interval $[a,x_1]\subset[a,b].$ Choose $x_1$ such that
\begin{equation}\label{a1}
  A\frac{\Gamma(\gamma)}{\Gamma(\alpha+\gamma)}(x_1-a)^{\alpha}<1
\end{equation}
holds and first we prove the existence of unique solution $y\in C_{1-\gamma}[a,x_1].$ We proceed as follows. Set Picard's sequence functions
\begin{equation}\label{a2}
  y_{0}(x)=\frac{y_a}{\Gamma(\gamma)}(x-a)^{\gamma-1},\qquad\gamma=\alpha+\beta-\alpha\beta,
\end{equation}
\begin{equation}\label{a3}
  y_m(x)=y_0(x)+\frac{1}{\Gamma(\alpha)}\int_{a}^{x}(x-t)^{\alpha-1}f(t,y_{m-1}(t))dt, \qquad m\in\N.
\end{equation}
We now show that $y_m(x)\in C_{1-\gamma}[a,b].$ From equation \eqref{a2}, it follows that $y_0(x)\in C_{1-\gamma}[a,b].$ By Lemma 3.1, $I_{a^{+}}^{\alpha}f$ is bounded from $C_{1-\gamma}[a,b]$ to $C_{1-\gamma}[a,b],$ which gives $y_m(x)\in C_{1-\gamma}[a,b],m\in\N$ and $x\in(a,b].$

By equations \eqref{a2} and \eqref{a3}, we have
\begin{equation*}
{\|y_1(x)-y_{0}(x)\|}_{C_{1-\gamma}[a,x_1]}={\|I_{a^{+}}^{\alpha}f(x,y_{0}(x))\|}_{C_{1-\gamma}[a,x_1]}
\end{equation*}
by using Lemma 3.1,
\begin{align}\label{a4}
 {\|y_1(x)-y_{0}(x)\|}_{C_{1-\gamma}[a,x_1]}&\leq M \frac{\Gamma(\gamma)}{\Gamma(\gamma+\alpha)}(x_1-a)^{\alpha}.
\end{align}
Further we obtain
\begin{equation}\label{a5}
  {\|y_{2}(x)-y_{1}(x)\|}_{C_{1-\gamma}[a,x_1]}\leq M \frac{\Gamma(\gamma)}{\Gamma(\gamma+\alpha)}(x_1-a)^{\alpha}
  \bigg( A \frac{\Gamma(\gamma)}{\Gamma(\gamma+\alpha)}(x_1-a)^{\alpha}\bigg).
\end{equation}
Continuing in this way $m$-times, we obtain
\begin{equation}\label{a6}
\hspace{-0.25cm}{\|y_{m}(x)-y_{m-1}(x)\|}_{C_{1-\gamma}[a,x_1]}\leq\frac{M\Gamma(\gamma)}{\Gamma(\gamma+\alpha)}(x_1-a)^{\alpha}
{\bigg(\frac{A(x_1-a)^{\alpha}\Gamma(\gamma)}{\Gamma(\gamma+\alpha)}\bigg)}^{m-1}
\end{equation}
By equation \eqref{a1}, we get
\begin{equation}\label{a7}
  {\|y_{m}(x)-y(x)\|}_{C_{1-\gamma}[a,x_1]}\to0,\qquad \text{as}\quad m\to+\infty.
\end{equation}
Again by Lemma 3.1, it follows that
\begin{equation*}
  {\|I_{a^+}^{\alpha}f(x,y_{m}(x))-I_{a^+}^{\alpha}f(x,y(x))\|}_{C_{1-\gamma}[a,x_1]}\leq A \frac{(x_1-a)^{\alpha}
  \Gamma(\gamma)}{\Gamma(\gamma+\alpha)} {\|y_{m}(x)-y(x)\|}_{C_{1-\gamma}[a,x_1]}
\end{equation*}
and hence by equation \eqref{a7},
\begin{equation}\label{a8}
  {\|I_{a^+}^{\alpha}f(x,y_{m}(x))-I_{a^+}^{\alpha}f(x,y(x))\|}_{C_{1-\gamma}[a,x_1]}\to0,\qquad \text{as}\quad m\to+\infty.
\end{equation}
From equations \eqref{a7} and \eqref{a8}, it follows that $y(x)$ is the solution of integral equation \eqref{s3} in $C_{1-\gamma}[a,x_1].$

Now to show that the solution $y(x)$ is unique, consider there exists two solutions $y(x)$ and $z(x)$ of the integral equation \eqref{s3} on $[a,x_1].$
Substituting them into \eqref{s3} and using Lemma 2.1 with Lipschits condition \eqref{d4}, we get
\begin{align}\label{a9}
  {\|y(x)-z(x)\|}_{C_{1-\gamma}[a,x_1]}&={\|I_{a^+}^{\alpha}f(x,y(x))-I_{a^+}^{\alpha}f(x,z(x))\|}_{C_{1-\gamma}[a,x_1]}\nonumber\\
  &\leq A \frac{(x_1-a)^{\alpha}\Gamma(\gamma)}{\Gamma(\gamma+\alpha)} {\|y(x)-z(x)\|}_{C_{1-\gamma}[a,x_1]}
\end{align}
This yields $A \frac{\Gamma(\gamma)}{\Gamma(\gamma+\alpha)}(x_1-a)^{\alpha}\geq1,$ which contradicts to condition \eqref{a1}.
Thus there exists $y(x)=y_1(x)\in C_{1-\gamma}[a,x_1]$ as a unique solution on $[a,x_1].$

Next, consider the interval $[x_1,x_2],$ where $x_2=x_1+h_1, h_1>0$ such that $x_2<b.$ Now the integral equation \eqref{s3} takes the form
\begin{align}\label{a10}
 y(x)=&\frac{y_a}{\Gamma(\gamma)}(x-a)^{\gamma-1}+\frac{1}{\Gamma(\alpha)}\int_{x_1}^{x}(x-t)^{\alpha-1}f(t,y(t))dt\nonumber\\
   &+\frac{1}{\Gamma(\alpha)}\int_{a}^{x_1}(x-t)^{\alpha-1}f(t,y(t))dt, \qquad x\in[x_1,x_2].
\end{align}
Since the function $y(x)$ is uniquely defined on $[a,x_1],$ the last integral is known function and therefore the integral equation \eqref{a10} can be written in the form
\begin{align}\label{a11}
 y(x)=y_{0}^{*}(x)+\frac{1}{\Gamma(\alpha)}\int_{x_1}^{x}(x-t)^{\alpha-1}f(t,y(t))dt, \quad x\in[x_1,x_2],
\end{align}
where
\begin{equation}\label{a12}
  y_{0}^{*}(x)=\frac{y_a}{\Gamma(\gamma)}(x-a)^{\gamma-1}+\frac{1}{\Gamma(\alpha)}\int_{a}^{x_1}(x-t)^{\alpha-1}f(t,y(t))dt
\end{equation}
is the known function. Using the same argument as above, we deduce that there exist a unique solution $y(x)=y_2(x)\in C_{1-\gamma}[x_1,x_2]$ on $[x_1,x_2].$ Taking interval $[x_2,x_3],$ where $x_3=x_2+h_2, h_2>0$ such that $x_3<b,$ and repeating the above process, we obtain a unique solution $y(x)\in C_{1-\gamma}[a,b]$ of integral equation \eqref{s3} such that
$y(x)=y_j(x)\in C_{1-\gamma}[x_{j-1},x_{j}]$, for $ j=1,2,..,l,$ and $a=x_0<x_2<..<x_l=b.$

Using differential equation \eqref{s1} and Lipschitz condition \eqref{d4}, we obtain
\begin{align}\label{a13}
  {\|D_{a^+}^{\alpha,\beta}y_{m}(x)-D_{a^+}^{\alpha,\beta}y(x)\|}_{C_{1-\gamma}[a,b]}&={\|f(x,y_{m}(x))-f(x,y(x))\|}_{C_{1-\gamma}[a,b]}\nonumber\\
&\leq A {\|y_{m}(x)-y(x)\|}_{C_{1-\gamma}[a,b]}.
\end{align}
Clearly, \eqref{a7} and \eqref{a13} implies that $D_{a^+}^{\alpha,\beta}y(x)\in C_{1-\gamma}[a,b]$ and result follows.
\begin{remark}
For the case $\beta=0,$ Theorem 3.1 reduces to the Cauchy-type problem for Riemann-Liouville fractional differential equation (\cite{cj}, Theorem 3.5) whereas for $\beta= 1$ in Theorem 3.1 yields the Cauchy-type problem for Caputo fractional differential equation (\cite{km}, Theorem 2).
\end{remark}
\section{Continuous Dependence}
In this section, first we study the continuous dependence of solution of Cauchy-type problem for Hilfer fractional differential equation using generalized Gronwall inequality as a handy tool. Consider the IVP \eqref{s1}-\eqref{s2}, for $0<\alpha<1,$ $0\leq\beta\leq1,$ $a\leq x<b,$ $(b\leq+\infty)$ and $f:[a,b)\times\R\to\R.$ To present dependence of solution on the order, let us consider the solutions of two IVPs with the neighbouring orders. We need the following lemma.
\begin{lem}\cite{ygd}
Suppose $\beta>0$, $a(t)$ is nonnegative function locally integrable on $0\leq t<T$ for some $(T\leq+\infty)$ and $g(t)$ is nonnegative, nondecreasing continuous function defined on $0\leq t<T,$ $g(t)\leq M$ (constant), and suppose $u(t)$ is nonnegative and locally integrable on $0\leq t<T$ with
\begin{equation*}
u(t)\leq a(t)+g(t)\int_{0}^{t}(t-s)^{\beta-1}u(s)ds
\end{equation*}
on this interval. Then
\begin{equation*}
u(t)\leq a(t)+\int_{0}^{t}\bigg[\sum_{n=1}^{\infty}\frac{(g(t)\Gamma(\beta))^{n}}{\Gamma(n\beta)}(t-s)^{n\beta-1}a(s)\bigg]ds,\quad 0\leq t<T.
\end{equation*}
\end{lem}
\begin{theorem}
Let $\alpha>0,\delta>0$ such that $0<\alpha-\delta<\alpha\leq1.$ Let $f$ is continuous function satisfying Lipschitz condition \eqref{d4} in $\R.$ For $a\leq x\leq h<b,$ assume that $y$ is the solution of IVP \eqref{s1}-\eqref{s2} and
$\hat{y}$ is the solution of IVP
\begin{equation}\label{k1}
D_{a^+}^{\alpha-\delta,\beta}\hat{y}(x)=f(x,\hat{y}),\qquad 0<\alpha<1,\, 0\leq\beta\leq1,
\end{equation}
\begin{equation}\label{k2}
I_{a^+}^{1-\gamma-\delta(\beta-1)}\hat{y}(x){\big|}_{x=a}={\hat{y}}_a,\qquad \gamma=\alpha+\beta(1-\alpha).
\end{equation}
Then for $a < x \leq h,$
\begin{equation*}
|\hat{y}(x)-{y}(x)|\leq{B(x)}+\int_{a}^{x}\bigg[\sum_{n=1}^{\infty}{\bigg(\frac{A}{\Gamma(\alpha)}\Gamma(\alpha-\delta)\bigg)}^{n}\frac{(x-t)^{n(\alpha-\delta)-1}}{\Gamma(n(\alpha-\delta))}B(t)\bigg]dt
\end{equation*}
hold, where
\begin{align}\label{k3}
B(x)&=\bigg|\frac{{\hat{y}}_{a}(x-a)^{\gamma+\delta(\beta-1)-1}}{\Gamma(\gamma+\delta(\beta-1))}-\frac{{y}_{a}(x-a)^{\gamma-1}}{\Gamma(\gamma)}\bigg|+\nonumber\\
&\|f\|\bigg|\frac{(x-a)^{\alpha-\delta}}{\Gamma(\alpha-\delta+1)}-\frac{(x-a)^{\alpha-\delta}}{(\alpha-\delta)\Gamma(\alpha)}\bigg|
+\|f\|\bigg|\frac{(x-a)^{\alpha-\delta}}{(\alpha-\delta)\Gamma(\alpha)}-\frac{(x-a)^{\alpha}}{\Gamma(\alpha+1)}\bigg|
\end{align}
and
\begin{equation*}
\|f\|=\max_{a\leq x\leq h}|f(x,y(x))|.
\end{equation*}
\end{theorem}
\textbf{Proof:} The equivalent integral solutions of IVPs \eqref{s1}-\eqref{s2} and \eqref{k1}-\eqref{k2} are
\begin{equation*}
y(x)=\frac{y_a}{\Gamma(\gamma)}(x-a)^{\gamma-1}+\frac{1}{\Gamma(\alpha)}\int_{a}^{x}(x-t)^{\alpha-1}f(t,y(t))dt
\end{equation*}
and
\begin{equation*}
\hat{y}(x)=\frac{{\hat{y}}_a}{\Gamma(\gamma+\delta(\beta-1))}(x-a)^{\gamma+\delta(\beta-1)-1}+\frac{1}{\Gamma(\alpha-\delta)}\int_{a}^{x}(x-t)^{\alpha-\delta-1}f(t,\hat{y}(t))dt,
\end{equation*}
respectively. It follows that
\begin{align*}
|\hat{y}(x)-{y}(x)|&=\bigg|\frac{{\hat{y}}_{a}}{\Gamma(\gamma+\delta(\beta-1))}(x-a)^{\gamma+\delta(\beta-1)-1}-\frac{{y}_{a}}{\Gamma(\gamma)}(x-a)^{\gamma-1}\\
&\hspace{1.00cm}+\int_{a}^{x}\frac{(x-t)^{\alpha-\delta-1}}{\Gamma(\alpha-\delta)}f(t,\hat{y}(t))dt-\int_{a}^{x}\frac{(x-t)^{\alpha-1}}{\Gamma(\alpha)}f(t,y(t))dt\bigg|\\
&\leq\bigg|\frac{{\hat{y}}_{a}}{\Gamma(\gamma+\delta(\beta-1))}(x-a)^{\gamma+\delta(\beta-1)-1}-\frac{{y}_{a}}{\Gamma(\gamma)}(x-a)^{\gamma-1}\bigg|\\
&\hspace{1.00cm}+\bigg|\int_{a}^{x}\bigg[\frac{(x-t)^{\alpha-\delta-1}}{\Gamma(\alpha-\delta)}-\frac{(x-t)^{\alpha-\delta-1}}{\Gamma(\alpha)}\bigg]|f(t,\hat{y}(t))|dt\bigg|\\
&\hspace{1.2cm}+\bigg|\frac{1}{\Gamma(\alpha)}\int_{a}^{x}(x-t)^{\alpha-\delta-1}|f(t,\hat{y}(t))-f(t,y(t))|dt\bigg|\\
&\hspace{1.2cm}+\bigg|\frac{1}{\Gamma(\alpha)}\int_{a}^{x}\big[(x-t)^{\alpha-\delta-1}-(x-t)^{\alpha-1}\big] |f(t,y(t))|dt\bigg|\\
|\hat{y}(x)-{y}(x)|&\leq B(x)+\frac{A}{\Gamma(\alpha)}\int_{a}^{x}(x-t)^{\alpha-\delta-1}|\hat{y}(t)-y(t)|dt,
\end{align*}
where $B(x)$ is defined by \eqref{k3}. Applying Lemma 4.1, we obtain
\begin{equation*}
|\hat{y}(x)-{y}(x)|\leq{B(x)}+\int_{a}^{x}\bigg[\sum_{n=1}^{\infty}{\bigg(\frac{A}{\Gamma(\alpha)}\Gamma(\alpha-\delta)\bigg)}^{n}\frac{(x-t)^{n(\alpha-\delta)-1}}{\Gamma(n(\alpha-\delta))}B(t)\bigg]dt,
\end{equation*}
the result is proved.

Next, we consider the fractional differential equation \eqref{s1} with the small change in the initial condition \eqref{s2},
\begin{equation}\label{c1}
I_{a^{+}}^{1-\gamma}y(x){|}_{x=a}=y_a+{\epsilon},\quad \gamma=\alpha+\beta(1-\alpha),
\end{equation}
where ${\epsilon}$ is arbitrary constant. We state and prove the result as follows:
\begin{theorem}
Suppose that assumptions of Theorem 3.1 hold. Suppose $y(x)$ and $\hat{y}(x)$ are solutions of IVPs \eqref{s1}-\eqref{s2} and \eqref{s1}-\eqref{c1}, respectively. Then
\begin{equation}\label{c2}
|y(x)-\hat{y}(x)|\leq|{\epsilon}|(x-a)^{\gamma-1}E_{\alpha,\gamma}(A(x-a)^{\alpha}), \quad x\in(a,b],
\end{equation}
holds, where $E_{\alpha,\gamma}(z)=\sum_{k=0}^{\infty}\frac{z^k}{\Gamma(k\alpha+\gamma)}$ is the Mittag-Leffler function, (see \cite{ss}).
\end{theorem}
\textbf{Proof:} In accordance with Theorem 3.1, we have $y(x)=\lim_{m\to\infty}y_m(x),$ with $y_{0}(x)$ and $y_{m}(x)$ are as defined in equations \eqref{a2} and \eqref{a3}, respectively. Clearly, we can write $\hat{y}(x)=\lim_{m\to\infty}\hat{y}_m(x),$ and
\begin{equation}\label{c3}
\hat{y}_{0}(x)=\frac{(y_a+{\epsilon})}{\Gamma(\gamma)}(x-a)^{\gamma-1},
\end{equation}
\begin{equation}\label{c4}
\hat{y}_{m}(x)=\hat{y}_{0}(x)+\frac{1}{\Gamma(\alpha)}\int_{a}^{x}(x-t)^{\alpha-1}f(t,\hat{y}_{m-1}(t))dt.
\end{equation}
From \eqref{a2} and \eqref{c3}, we have
\begin{align}\label{c5}
\big|{y}_{0}(x)-\hat{y}_{0}(x)\big|&=\bigg|\frac{y_a}{\Gamma(\gamma)}(x-a)^{\gamma-1}-\frac{(y_a+{\epsilon})}{\Gamma(\gamma)}(x-a)^{\gamma-1}\bigg|\nonumber\\
\big|{y}_{0}(x)-\hat{y}_{0}(x)\big|&\leq  \big|{\epsilon}\big|\frac{(x-a)^{\gamma-1}}{\Gamma(\gamma)}.
\end{align}
By the subsequent relations \eqref{a3} and \eqref{c4}, the Lipschitz condition \eqref{d4} and the inequality \eqref{c5}, we obtain
\begin{align}\label{c6}
\big|{y}_{1}(x)-\hat{y}_{1}(x)\big|&=\big|{\epsilon}\frac{(x-a)^{\gamma-1}}{\Gamma(\gamma)}+\int_{a}^{x}
\frac{(x-t)^{\alpha-1}}{\Gamma(\alpha)}[f(t,y_{0}(t))-f(t,\hat{y}_{0}(t))] dt\big| \nonumber\\
&\leq|{\epsilon}|\frac{(x-a)^{\gamma-1}}{\Gamma(\gamma)}+\frac{A}{\Gamma(\alpha)}\int_{a}^{x}(x-t)^{\alpha-1}
|y_{0}(t)-\hat{y}_{0}(t)|dt\nonumber\\
&\leq |{\epsilon}|\frac{(x-a)^{\gamma-1}}{\Gamma(\gamma)}+A|{\epsilon}|
\frac{(x-a)^{\alpha+\gamma-1}}{\Gamma(\alpha+\gamma)}\nonumber\\
\big|{y}_{1}(x)-\hat{y}_{1}(x)\big|&\leq|{\epsilon}|(x-a)^{\gamma-1}\sum_{j=0}^{1}\frac{A^j(x-a)^{\alpha{j}}}{\Gamma(\alpha{j}+\gamma)}.
\end{align}
Similarly, by using \eqref{c6}, it directly follows that
\begin{align*}
\big|{y}_{2}(x)-\hat{y}_{2}(x)\big|&\leq|{\epsilon}|\frac{(x-a)^{\gamma-1}}{\Gamma(\gamma)}+\frac{A}{\Gamma(\alpha)}\int_{a}^{x}(x-t)^{\alpha-1}
|y_{1}(t)-\hat{y}_{1}(t)|dt\\
&\leq|{\epsilon}|(x-a)^{\gamma-1}\sum_{j=0}^{2}\frac{A^j(x-a)^{\alpha{j}}}{\Gamma(\alpha{j}+\gamma)}.
\end{align*}
By the induction, we obtain
\begin{equation}\label{c7}
\big|{y}_{m}(x)-\hat{y}_{m}(x)\big|\leq|{\epsilon}|(x-a)^{\gamma-1}\sum_{j=0}^{m}\frac{A^j(x-a)^{\alpha{j}}}{\Gamma(\alpha{j}+\gamma)}.
\end{equation}
Taking limit as $m\to\infty$ in \eqref{c7}, we have
\begin{align*}
\big|y(x)-\hat{y}(x)\big|&\leq|{\epsilon}|(x-a)^{\gamma-1}\sum_{j=0}^{\infty}\frac{A^j(x-a)^{\alpha{j}}}{\Gamma(\alpha{j}+\gamma)}\\
&=|{\epsilon}|(x-a)^{\gamma-1}E_{\alpha,\gamma}(A(x-a)^{\alpha}),
\end{align*}
which completes the proof of Theorem 4.2.
\begin{remark}
It follows from Theorem 4.1 and 4.2 that, small change in order and initial condition \eqref{s2} cause only small change in the solution on $[l,b]$  for $l$ between $a$ to $b$ which does not contain initial point $a.$ On the other hand, the solution may change significantly in the interval $[a,l].$
\end{remark}

\end{document}